\documentstyle[12pt,leqno]{article}

\catcode`\@=11
\def\relaxnext@{\let\next\relax}
\font\tenmsx=msxm10 scaled\magstep1
\font\sevenmsx=msxm7 scaled\magstep1
\font\fivemsx=msxm5 scaled\magstep1
\font\tenmsy=msym10 scaled\magstep1
\font\sevenmsy=msym7 scaled\magstep1
\font\fivemsy=msym5 scaled\magstep1
\newfam\msxfam
\newfam\msyfam
\textfont\msxfam=\tenmsx  \scriptfont\msxfam=\sevenmsx
  \scriptscriptfont\msxfam=\fivemsx
\textfont\msyfam=\tenmsy  \scriptfont\msyfam=\sevenmsy
  \scriptscriptfont\msyfam=\fivemsy

\def\hexnumber@#1{\ifcase#1 0\or1\or2\or3\or4\or5\or6\or7\or8\or9\or
	A\or B\or C\or D\or E\or F\fi }

\font\teneuf=eufm10 scaled\magstep1
\font\seveneuf=eufm7 scaled\magstep1
\font\fiveeuf=eufm5 scaled\magstep1
\newfam\euffam
\textfont\euffam=\teneuf
\scriptfont\euffam=\seveneuf
\scriptscriptfont\euffam=\fiveeuf
\def\frak{\relaxnext@\ifmmode\let\next\frak@\else
 \def\next{\errmessage{Use \string\frak\space only in math
mode}}\fi\next}
\def\goth{\relaxnext@\ifmmode\let\next\frak@\else 
\def\next{\errmessage{Use \string\goth\space only in math
mode}}\fi\next} 
\def\frak@#1{{\frak@@{#1}}}
\def\frak@@#1{\fam\euffam#1}

\edef\msx@{\hexnumber@\msxfam}
\edef\msy@{\hexnumber@\msyfam}

\mathchardef\boxdot="2\msx@00
\mathchardef\boxplus="2\msx@01
\mathchardef\boxtimes="2\msx@02
\mathchardef\square="0\msx@03
\mathchardef\blacksquare="0\msx@04
\mathchardef\centerdot="2\msx@05
\mathchardef\lozenge="0\msx@06
\mathchardef\blacklozenge="0\msx@07
\mathchardef\circlearrowright="3\msx@08
\mathchardef\circlearrowleft="3\msx@09
\mathchardef\rightleftharpoons="3\msx@0A
\mathchardef\leftrightharpoons="3\msx@0B
\mathchardef\boxminus="2\msx@0C
\mathchardef\Vdash="3\msx@0D
\mathchardef\Vvdash="3\msx@0E
\mathchardef\vDash="3\msx@0F
\mathchardef\twoheadrightarrow="3\msx@10
\mathchardef\twoheadleftarrow="3\msx@11
\mathchardef\leftleftarrows="3\msx@12
\mathchardef\rightrightarrows="3\msx@13
\mathchardef\upuparrows="3\msx@14
\mathchardef\downdownarrows="3\msx@15
\mathchardef\upharpoonright="3\msx@16

\mathchardef\downharpoonright="3\msx@17
\mathchardef\upharpoonleft="3\msx@18
\mathchardef\downharpoonleft="3\msx@19
\mathchardef\rightarrowtail="3\msx@1A
\mathchardef\leftarrowtail="3\msx@1B
\mathchardef\leftrightarrows="3\msx@1C
\mathchardef\rightleftarrows="3\msx@1D
\mathchardef\Lsh="3\msx@1E
\mathchardef\Rsh="3\msx@1F
\mathchardef\rightsquigarrow="3\msx@20
\mathchardef\leftrightsquigarrow="3\msx@21
\mathchardef\looparrowleft="3\msx@22
\mathchardef\looparrowright="3\msx@23
\mathchardef\circeq="3\msx@24
\mathchardef\succsim="3\msx@25
\mathchardef\gtrsim="3\msx@26
\mathchardef\gtrapprox="3\msx@27
\mathchardef\multimap="3\msx@28
\mathchardef\therefore="3\msx@29
\mathchardef\because="3\msx@2A
\mathchardef\doteqdot="3\msx@2B

\mathchardef\triangleq="3\msx@2C
\mathchardef\precsim="3\msx@2D
\mathchardef\lesssim="3\msx@2E
\mathchardef\lessapprox="3\msx@2F
\mathchardef\eqslantless="3\msx@30
\mathchardef\eqslantgtr="3\msx@31
\mathchardef\curlyeqprec="3\msx@32
\mathchardef\curlyeqsucc="3\msx@33
\mathchardef\preccurlyeq="3\msx@34
\mathchardef\leqq="3\msx@35
\mathchardef\leqslant="3\msx@36
\mathchardef\lessgtr="3\msx@37
\mathchardef\backprime="0\msx@38
\mathchardef\risingdotseq="3\msx@3A
\mathchardef\fallingdotseq="3\msx@3B
\mathchardef\succcurlyeq="3\msx@3C
\mathchardef\geqq="3\msx@3D
\mathchardef\geqslant="3\msx@3E
\mathchardef\gtrless="3\msx@3F
\mathchardef\sqsubset="3\msx@40
\mathchardef\sqsupset="3\msx@41
\mathchardef\vartriangleright="3\msx@42
\mathchardef\vartriangleleft="3\msx@43
\mathchardef\trianglerighteq="3\msx@44
\mathchardef\trianglelefteq="3\msx@45
\mathchardef\bigstar="0\msx@46
\mathchardef\between="3\msx@47
\mathchardef\blacktriangledown="0\msx@48
\mathchardef\blacktriangleright="3\msx@49
\mathchardef\blacktriangleleft="3\msx@4A
\mathchardef\vartriangle="0\msx@4D
\mathchardef\blacktriangle="0\msx@4E
\mathchardef\triangledown="0\msx@4F
\mathchardef\eqcirc="3\msx@50
\mathchardef\lesseqgtr="3\msx@51
\mathchardef\gtreqless="3\msx@52
\mathchardef\lesseqqgtr="3\msx@53
\mathchardef\gtreqqless="3\msx@54
\mathchardef\Rrightarrow="3\msx@56
\mathchardef\Lleftarrow="3\msx@57
\mathchardef\veebar="2\msx@59
\mathchardef\barwedge="2\msx@5A
\mathchardef\doublebarwedge="2\msx@5B
\mathchardef\angle="0\msx@5C
\mathchardef\measuredangle="0\msx@5D
\mathchardef\sphericalangle="0\msx@5E
\mathchardef\varpropto="3\msx@5F
\mathchardef\smallsmile="3\msx@60
\mathchardef\smallfrown="3\msx@61
\mathchardef\Subset="3\msx@62
\mathchardef\Supset="3\msx@63
\mathchardef\Cup="2\msx@64

\mathchardef\Cap="2\msx@65

\mathchardef\curlywedge="2\msx@66
\mathchardef\curlyvee="2\msx@67
\mathchardef\leftthreetimes="2\msx@68
\mathchardef\rightthreetimes="2\msx@69
\mathchardef\subseteqq="3\msx@6A
\mathchardef\supseteqq="3\msx@6B
\mathchardef\bumpeq="3\msx@6C
\mathchardef\Bumpeq="3\msx@6D
\mathchardef\lll="3\msx@6E

\mathchardef\ggg="3\msx@6F

\mathchardef\circledS="0\msx@73
\mathchardef\pitchfork="3\msx@74
\mathchardef\dotplus="2\msx@75
\mathchardef\backsim="3\msx@76
\mathchardef\backsimeq="3\msx@77
\mathchardef\complement="0\msx@7B
\mathchardef\intercal="2\msx@7C
\mathchardef\circledcirc="2\msx@7D
\mathchardef\circledast="2\msx@7E
\mathchardef\circleddash="2\msx@7F
\def\ulcorner{\delimiter"4\msx@70\msx@70 }
\def\urcorner{\delimiter"5\msx@71\msx@71 }
\def\llcorner{\delimiter"4\msx@78\msx@78 }
\def\lrcorner{\delimiter"5\msx@79\msx@79 }
\def\yen{\mathhexbox\msx@55 }
\def\checkmark{\mathhexbox\msx@58 }
\def\circledR{\mathhexbox\msx@72 }
\def\maltese{\mathhexbox\msx@7A }
\mathchardef\lvertneqq="3\msy@00
\mathchardef\gvertneqq="3\msy@01
\mathchardef\nleq="3\msy@02
\mathchardef\ngeq="3\msy@03
\mathchardef\nless="3\msy@04
\mathchardef\ngtr="3\msy@05
\mathchardef\nprec="3\msy@06
\mathchardef\nsucc="3\msy@07
\mathchardef\lneqq="3\msy@08
\mathchardef\gneqq="3\msy@09
\mathchardef\nleqslant="3\msy@0A
\mathchardef\ngeqslant="3\msy@0B
\mathchardef\lneq="3\msy@0C
\mathchardef\gneq="3\msy@0D
\mathchardef\npreceq="3\msy@0E
\mathchardef\nsucceq="3\msy@0F
\mathchardef\precnsim="3\msy@10
\mathchardef\succnsim="3\msy@11
\mathchardef\lnsim="3\msy@12
\mathchardef\gnsim="3\msy@13
\mathchardef\nleqq="3\msy@14
\mathchardef\ngeqq="3\msy@15
\mathchardef\precneqq="3\msy@16
\mathchardef\succneqq="3\msy@17
\mathchardef\precnapprox="3\msy@18
\mathchardef\succnapprox="3\msy@19
\mathchardef\lnapprox="3\msy@1A
\mathchardef\gnapprox="3\msy@1B
\mathchardef\nsim="3\msy@1C
\mathchardef\ncong="3\msy@1D

\mathchardef\varsubsetneq="3\msy@20
\mathchardef\varsupsetneq="3\msy@21
\mathchardef\nsubseteqq="3\msy@22
\mathchardef\nsupseteqq="3\msy@23
\mathchardef\subsetneqq="3\msy@24
\mathchardef\supsetneqq="3\msy@25
\mathchardef\varsubsetneqq="3\msy@26
\mathchardef\varsupsetneqq="3\msy@27
\mathchardef\subsetneq="3\msy@28
\mathchardef\supsetneq="3\msy@29
\mathchardef\nsubseteq="3\msy@2A
\mathchardef\nsupseteq="3\msy@2B
\mathchardef\nparallel="3\msy@2C
\mathchardef\nmid="3\msy@2D
\mathchardef\nshortmid="3\msy@2E
\mathchardef\nshortparallel="3\msy@2F
\mathchardef\nvdash="3\msy@30
\mathchardef\nVdash="3\msy@31
\mathchardef\nvDash="3\msy@32
\mathchardef\nVDash="3\msy@33
\mathchardef\ntrianglerighteq="3\msy@34
\mathchardef\ntrianglelefteq="3\msy@35
\mathchardef\ntriangleleft="3\msy@36
\mathchardef\ntriangleright="3\msy@37
\mathchardef\nleftarrow="3\msy@38
\mathchardef\nrightarrow="3\msy@39
\mathchardef\nLeftarrow="3\msy@3A
\mathchardef\nRightarrow="3\msy@3B
\mathchardef\nLeftrightarrow="3\msy@3C
\mathchardef\nleftrightarrow="3\msy@3D
\mathchardef\divideontimes="2\msy@3E
\mathchardef\varnothing="0\msy@3F
\mathchardef\nexists="0\msy@40
\mathchardef\mho="0\msy@66
\mathchardef\eth="0\msy@67
\mathchardef\eqsim="3\msy@68
\mathchardef\beth="0\msy@69
\mathchardef\gimel="0\msy@6A
\mathchardef\daleth="0\msy@6B
\mathchardef\lessdot="3\msy@6C
\mathchardef\gtrdot="3\msy@6D
\mathchardef\ltimes="2\msy@6E
\mathchardef\rtimes="2\msy@6F
\mathchardef\shortmid="3\msy@70
\mathchardef\shortparallel="3\msy@71
\mathchardef\smallsetminus="2\msy@72
\mathchardef\thicksim="3\msy@73
\mathchardef\thickapprox="3\msy@74
\mathchardef\approxeq="3\msy@75
\mathchardef\succapprox="3\msy@76
\mathchardef\precapprox="3\msy@77
\mathchardef\curvearrowleft="3\msy@78
\mathchardef\curvearrowright="3\msy@79
\mathchardef\digamma="0\msy@7A
\mathchardef\varkappa="0\msy@7B
\mathchardef\hslash="0\msy@7D
\mathchardef\hbar="0\msy@7E
\mathchardef\backepsilon="3\msy@7F
\def\Bbb{\ifmmode\let\next\Bbb@\else
 \def\next{\errmessage{Use \string\Bbb\space only in math mode}}\fi\next}
\def\Bbb@#1{{\Bbb@@{#1}}}
\def\Bbb@@#1{\fam\msyfam#1}

\catcode`\@=12

\newcommand{\spline}{\vspace{1mm}}
\newcommand{\si}{\sigma}

\begin{document}

\title{KAEHLER STRUCTURES ON $K_{\bf C}/(P,P)$}

\author{MENG-KIAT CHUAH}

\date{}

\maketitle

\begin{abstract}
Let $K$ be a compact connected semi-simple Lie group, let $G = K_{\bf C}$, 
and let $G = KAN$ be an Iwasawa decomposition. 
Given a $K$-invariant Kaehler
structure $\omega$ on $G/N$, there corresponds a pre-quantum line bundle
{\bf L} on $G/N$. 
Following a suggestion of A.S.~Schwarz,
in a joint work with V. Guillemin, 
we study
its holomorphic sections ${\cal O}({\bf L})$ as a 
$K$-representation space.
We define a $K$-invariant $L^2$-structure on ${\cal O}({\bf L})$,
and let $H_\omega \subset {\cal O}({\bf L})$ denote the space
of square-integrable holomorphic sections. Then $H_\omega$ is a unitary
$K$-representation space, 
but we find
that not all unitary irreducible $K$-representations occur
as subrepresentations of $H_\omega$.
This paper serves as a continuation of that work,
by generalizing the space considered. Instead of working with $G/N = G/(B,B)$,
where $B$ is a Borel subgroup containing $N$,
we consider $G/(P,P)$, for all parabolic subgroups $P$ containing $B$.
We carry out similar construction, and recover 
in $H_\omega$ the unitary irreducible
$K$-representations previously missing. As a result, we use these
holomorphic sections to construct a model for $K$: a unitary $K$-representation
in which every irreducible $K$-representation occurs with multiplicity one.
\end{abstract}


\begin{center}
\section{INTRODUCTION}
\end{center}
\setcounter{equation}{0}

Let $K$ be a compact connected semi-simple Lie group, let $G = K_{\bf C}$
be its complexification, and let $G = KAN$ be an Iwasawa decomposition.
Since $G$ and $N$ are complex Lie groups, $G/N$ is a complex manifold,
and $G$ acts on $G/N$ by left action. Let $T$ be the centralizer of $A$ in
$K$, so that $H = TA$ is a Cartan subgroup of $G$. Since $H$ normalizes
$N$, there is a right action of $H$ on $G/N$. We shall often be interested
in the maximal compact group action of $K \times T$. We let
$\frak g, \frak k, \frak h,\frak t,\frak a, \frak n$ 
denote the Lie algebras of
$G, K, H, T, A, N$ respectively.

The following scheme of geometric quantization was suggested 
by A.S.~Schwarz  \cite{kn:lns}:
Equip $G/N$ with a suitable $K$-invariant Kaehler structure $\omega$, and
consider the pre-quantum line bundle ${\bf L}$ associated to
$\omega$ (\cite{kn:gs}, \cite{kn:ko}).
The Chern class of ${\bf L}$ is
$[\omega]$, and ${\bf L}$ comes with a connection $\nabla$ whose curvature
is $\omega$, as well as an invariant Hermitian structure $<,>$.   
We denote by ${\cal O}({\bf L})$ the space of holomorphic sections on
${\bf L}$. The $K$-action on $G/N$ lifts to a $K$-representation on
${\cal O}({\bf L})$.
Let $\mu$ be the $K \times A$-invariant measure on $G/N$, which is unique
up to non-zero constant. Given a holomorphic section $s$ of ${\bf L}$,
we consider the integral
 \[ \int_{G/N} <s,s> \mu \; .\]
Let $H_\omega \subset {\cal O}({\bf L})$ denote the holomorphic sections
in which this integral converges. Since $\mu$ is $K$-invariant,
$H_\omega$ becomes a unitary $K$-representation space.
It was hoped in \cite{kn:lns} that every irreducible
$K$-representation occurs with multiplicity one in $H_\omega$
(called a {\em model} by I.M.~Gelfand~\cite{kn:gz}). 

By the method of highest weight, the irreducible
$K$-representations 
can be labeled by the dominant integral weights in $\frak t^*$,
up to isomorphism. 
In a joint work with V.~Guillemin \cite{kn:cg}, 
we carry out this construction, but
find that no matter how $\omega$ is chosen,
the irreducibles whose highest weights
lie on the wall of the Weyl chamber do not occur in the Hilbert space
$H_\omega$.
Therefore, not all
unitary $K$-irreducibles occur in $H_\omega$.
The present paper follows a suggestion of V.~Guillemin (\cite{kn:cg} p.192),
by modifying the space $G/N$ to more general classes of
homogeneous spaces.
As a result, we manage to recover the unitary $K$-irreducibles
previously missing.  

Let $B = HN$ be the Borel subgroup of $G$. Observe that $(B,B) = N$, hence
$G/N = G/(B,B)$. With this in mind, 
we can generalize the class of homogeneous spaces considered to $G/(P,P)$,
 for $P$ a parabolic subgroup of $G$ containing $B$.
Since $P$ is a complex Lie group, so is $(P,P)$; hence $G/(P,P)$ is a
complex manifold. Clearly $G$ acts on $G/(P,P)$ on the left, 
and we shall see that 
a complex subgroup of $H$ normalizes $(P,P)$, and hence acts 
on $G/(P,P)$ on the right.

Let $W \subset \frak t^*$ denote the open Weyl chamber, and 
$\overline{W}$ its closure. We say that $\si \subset \overline{W}$ is a 
{\em cell} if there exists a subset $S$ of the positive simple roots $\Delta$
such that
\begin{equation}
 \si = \{x \in \overline{W} \;;\; 
(x,S)=0 \,,\, (x,\Delta \backslash S)>0 \} , 
\label{eq:cell}
\end{equation}
where the pairing used is the Killing form.
This way, $\overline{W}$ is a disjoint union of the cells
of various dimensions.
Using the Killing form and 
the almost complex structure,
it is convenient to regard the cell $\si$ as contained in any of the spaces
$\frak h, \frak t, \frak a, \frak h^*, \frak t^*, \frak a^*$,
depending on the context.
The cell $\si$ defines a subalgebra $\frak h_\si$ of $\frak h$, by taking
complex linear span of $\si$. 
Similarly, the subalgebras $\frak t_\si, \frak a_\si$
are defined by intersecting $\frak h_\si$ with $\frak t, \frak a$ respectively.
These subalgebras define the subgroups $H_\si, T_\si, A_\si$ of
$H, T, A$ respectively. 
A bijective correspondence between the cells 
$\{\si\}$ and the parabolic subgroups $\{P\}$
containing $B$ is given by Langlands decomposition (\cite{kn:kn} p.132)
\begin{equation}
 P = M A_\si N_\si . 
\label{eq:man}
\end{equation}

Fix a parabolic subgroup $P$ containing $B$,
with $\si$ its corresponding cell. 
Since $H_\si$ is the normalizer of $(P,P)$ in $H$, it acts on $G/(P,P)$ 
on the right. Out of the action of the complex group $G \times H_\si$, 
we shall consider the action of the maximal compact group $K \times T_\si$ on
$G/(P,P)$. We shall show that

$\spline$

\noindent {\bf Theorem I  } {\em  Let } $\omega$ {\em be a K-invariant Kaehler
structure on }
$G/(P,P).$ {\em  Then } $\omega$ {\em is } $K \times T_\si$
{\em -invariant if and only if it has a potential function.}
 
$\spline$

Though we shall be interested mostly in Kaehler structures, 
Theorem I holds also
for degenerate {\em (1,1)}-form $\omega$. 
In the next theorem, we shall derive a necessary and sufficient condition
for a {\em (1,1)}-form $\omega$ to be Kaehler.
Let $\omega$ be a 
$K \times T_\si$-invariant {\em (1,1)}-form, so that
\[ \omega = \sqrt{-1} \partial \bar{\partial} F ,\]
for some function $F$ on $G/(P,P)$. Averaging by the compact group $K$
if necessary, we may assume that $F$ is $K$-invariant.
Let $K^\si$ be the centralizer of $T_\si$ in $K$.
It defines a compact semi-simple subgroup $K_{ss}^\si$ of $K$,
given by $K_{ss}^\si = (K^\si,K^\si)$.
We shall show that, as real manifolds and $K \times H_\si$-spaces,
\begin{equation}
 G/(P,P) = (K/K_{ss}^\si) \times A_\si .
\label{eq:coca}
\end{equation}
Therefore, the potential function
 $F$, being $K$-invariant, can be regarded as a function on $A_\si$.
Since the exponential map identifies the vector space $\frak a_\si$ with
$A_\si$, $F$ becomes a function on $\frak a_\si$. The almost complex structure
identifies the dual spaces $\frak a_\si^* \cong \frak t_\si^*$, hence the
Legendre transform of $F$ can be written as
\[ L_F : \frak a_\si \longrightarrow \frak t_\si^* .\]
The significance of this map will become apparent shortly, when we study
the moment map. We write $\log : A_\si \longrightarrow \frak a_\si$ for
the inverse of the exponential map.

The $K$-action on $G/(P,P)$ preserving $\omega$
is Hamiltonian: there exists a unique moment map
\[ \Phi : G/(P,P) \longrightarrow \frak k^* \]
corresponding to this action. Since $\Phi$ is $K$-equivariant, 
(\ref{eq:coca}) implies that it is
determined by its value on $A_\si \subset (K/K_{ss}^\si)\times A_\si$,
where $A_\si$ is imbedded as its product with the identity coset
$eK_{ss}^\si \in K/K_{ss}^\si$.
Meanwhile, since $\frak k$ is semi-simple, the Killing form on $\frak k$
is non-degenerate;
which induces the inclusion $\frak t^* \subset \frak k^*$ from
$\frak t \subset \frak k$.  

$\spline$

\noindent {\bf Theorem II  }
{\em  Let } $\omega$ {\em be a }
$K \times T_\si$ {\em -invariant (1,1)-form on } $G/(P,P)$.
{\em  Then its moment map } $\Phi$
{\em and its potential function F satisfy } 
$\Phi(a) = \frac{1}{2} L_F (\log a) \in \frak t_\si^*$ 
{\em  for all } $a \in A_\si$. 
{\em  Further, } $ \omega = \sqrt{-1}\partial \bar{\partial}F $
{\em  is Kaehler if and only if: }

{\em (i) }$\;\;$ $F \in C^\infty(\frak a_\si)$ {\em  is strictly convex; and }

{\em (ii)}$\;$ {\em  The image of } $ \frac{1}{2} L_F $
{\em  is contained in the cell } $ \si \subset \frak t_\si^*$;
{\em i.e.} $\Phi(A_\si) \subset \si$.

$\spline$

Since a $K \times T_\si$-invariant Kaehler structure $\omega$ has a
potential function $F$, it is exact.
Therefore, it is in particular integral.
Let ${\bf L}$ be the line bundle on $G/(P,P)$ whose Chern class is
$[\omega]=0$, 
equipped with a connection $\nabla$ whose curvature is $\omega$
(\cite{kn:gs},\cite{kn:ko}).
The topology of ${\bf L}$ is trivial, but the connection $\nabla$
gives rise to interesting geometry on the holomorphic sections of ${\bf L}$.
We recall that ${\bf L}$ is equipped with an invariant Hermitian
structure
$<,>$. Let $\mu$
be a $K \times A_\si$-invariant measure on $G/(P,P)$. 
We consider the integral 
\begin{equation}
\int_{G/(P,P)} <s,s> \mu \; , 
\label{eq:int}
\end{equation}
 for holomorphic
sections $s$ of ${\bf L}$. 
As we shall see in Theorem III, convergence of this integral 
is determined by the image of the moment map. The 
$K \times T_\si$-action on $G/(P,P)$
lifts to a $K \times T_\si$-representation on ${\cal O}({\bf L})$, the space of
holomorphic sections of ${\bf L}$. 
We similarly define $H_\omega \subset {\cal O}({\bf L})$ to be the 
holomorphic sections in which (\ref{eq:int}) converges.
Since $\mu$ is $K$-invariant, $H_\omega$ becomes a 
unitary $K$-representation space.
For a dominant integral weight $\lambda$,
let ${\cal O}({\bf L})_\lambda$ be the holomorphic
sections in ${\bf L}$ that transform by $\lambda$ under 
the right $T_\si$-action.
Since the left $K$-action commutes with the right $T_\si$-action,
${\cal O}({\bf L})_\lambda$ is a $K$-representation space.
Let $\si$ be the cell corresponding to the parabolic
subgroup $P$, and let $\overline{\si}$ be its closure. Then

$\spline$

\noindent {\bf Theorem III  } 
{\em  The irreducible K-representation with highest weight } $\lambda$
{\em occurs in } ${\cal O}({\bf L})$
{\em if and only if } $\lambda \in \overline{\si}$.
{\em For } $\lambda \in \overline{\si}$,
{\em it occurs with multiplicity one, and is given by }
${\cal O}({\bf L})_\lambda$.
{\em Further, } ${\cal O}({\bf L})_\lambda$ {\em is contained in } $H_\omega$
{\em if and only if } $\lambda$
{\em lies in the image of the moment map. }

$\spline$

With this result, it is now clear that in \cite{kn:cg}, the singular 
representations are never contained in $H_\omega$ :

When $P=B$, $\si$ 
becomes the open Weyl chamber $W$. 
Then Theorem II says that $\Phi(A_\si) \subset W$; 
and by $K$-equivariance, 
$\Phi(G/(P,P))= Ad_K^*(\Phi(A_\si))$
does not intersect the wall $\overline{W} \backslash W$.
Consequently, by Theorem III, the irreducible representations
${\cal O}({\bf L})_\lambda$ with highest
weight $\lambda \in \overline{W} \backslash W$ cannot be
contained in $H_\omega$.
 
Similarly, for general parabolic subgroup $P$, not all 
${\cal O}({\bf L})_\lambda$ are contained in $H_\omega$:
For $\lambda \in \overline{\si} \backslash \si$,
Theorems II and III say that ${\cal O}({\bf L})_\lambda$ exists non-trivially
but is not contained in $H_\omega$.

We shall see that, however,
for a suitable Kaehler structure $\omega$ on $G/(P,P)$, the image of the
moment map intersects $\overline{\si}$ in all of $\si$. This way, 
by Theorem III, all
the $K$-irreducibles ${\cal O}({\bf L})_\lambda$
with highest weights $\lambda \in \si$ are contained in $H_\omega$.
As an application,
we provide a geometric construction of a unitary
$K$-representation, containing all the irreducibles with multiplicity one.

$\spline$

\noindent {\bf Acknowledgement } The author would like to thank 
V. Guillemin, R. Sjamaar and D. Vogan for many helpful suggestions.
The referee has helped to clarify some definitions and notations used in
this paper.

\newpage
\begin{center}
\section{KAEHLER STRUCTURES ON $G/(P,P)$}
\end{center}
\setcounter{equation}{0}

The main purpose of this section is to prove Theorem I.
Since $K$ is connected and semi-simple, so is $G = K_{\bf C}$.
Let $P$ be a parabolic subgroup of $G$ containing $B$, 
and $\si$ the cell corresponding
to $P$. They are related by Langlands decomposition (\ref{eq:man})
\[ P = M A_\si N_\si \,, \]
where $A_\si$ is the subgroup described in $\S 1$.
 Then $A_\si \subset A, N_\si \subset N$, where $A, N$ come from Iwasawa
decomposition of $G$.
Further, $A_\si$ normalizes $N_\si$, and is the centralizer of 
$M A_\si$ in $A$. Therefore, $H_\si = T_\si A_\si$ is the normalizer of
$(P,P) = (M,M)N_\si$ in $H$, which induces a natural right $H_\si$-action
on $G/(P,P)$. We shall give another description of $G/(P,P)$, which
reflects this right action better.

Since $G$ is semi-simple, the Killing form is non-degenerate.
Let $\frak a_\si^\perp$ be the orthocomplement of $\frak a_\si$ 
with respect to the Killing form in $\frak a$,
and $A_\si^\perp \subset A$ the corresponding subgroup
induced by $\frak a_\si^\perp$. We construct 
$\frak t_\si^\perp, T_\si^\perp, \frak h_\si^\perp, H_\si^\perp$ similarly.
Let $K^\si$ be the subgroup of $K$ given by
\[ K^\si = \{ k \in K \;;\; kt = tk \mbox{ for all } t \in T_\si \} . \]
Let $K_{ss}^\si = (K^\si, K^\si)$ be the corresponding compact semi-simple
Lie group. Then
\begin{equation}
 (K_{ss}^\si)_{\bf C} = K_{ss}^\si \, A_\si^\perp \, (M \cap N) 
\label{eq:iwa}
\end{equation}
is an Iwasawa decomposition of the complexified group $(K_{ss}^\si)_{\bf C}$.
Since $N~=~(M~\cap~N)~N_\si$, it follows from (\ref{eq:iwa}) that
\begin{equation}
\begin{array}{rl}
K_{ss}^\si A_\si^\perp N & = (K_{ss}^\si)_{\bf C} N_\si \\
& = (K_{\bf C}^\si)_{ss} N_\si \\
& = (M A_\si , M A_\si) N_\si \\
& = (M, M) N_\si \\
& = (P,P) .
\end{array}
\label{eq:reyer}
\end{equation}
Then, the Iwasawa decomposition $G = KAN$ and (\ref{eq:reyer}) imply that
\begin{equation}
G/(P,P) = (K/K_{ss}^\si)\times A_\si ,
\label{eq:nnsc}
\end{equation}
as real manifolds and $K \times H_\si$-spaces.
With this description, the right action of $H_\si = T_\si A_\si$ is clear:
$T_\si$ acts on $(K/K_{ss}^\si)\times A_\si$
simply because it commutes with $K_{ss}^\si$ and $A_\si$,
while $A_\si$ acts on
$(K/K_{ss}^\si)\times A_\si$ by group multiplication on itself. We shall be concerned with the
$K \times T_\si$-action on $G/(P,P)$.

Since $N = (B,B) \subset (P,P)$, there is a fibration
\begin{equation}
\pi : G/N \longrightarrow G/(P,P) .
\label{eq:fib}
\end{equation}
It follows from $G=KAN$ and
(\ref{eq:nnsc}) that the fiber of $\pi$ is 
$K_{ss}^\si \times A_\si^\perp$. Further, $\pi$ sends every right $H$-orbit
in $G/N$ to a right $H_\si$-orbit in $G/(P,P)$, by contracting each
$H_\si^\perp$-coset to a point.
   
Given a $K$-invariant
Kaehler structure $\omega$ on $G/(P,P)$, 
we want to show that it is invariant under the right $T_\si$-action
if and only if it has a potential function.
Our strategy is to work on
the {\em (1,1)}-form $\pi^* \omega$ on $G/N$ using results in \cite{kn:cg},
then transfer this result back to $\omega$. 
 Let $V$ be the orthocomplement of $\frak t$ in $\frak k$
with respect to the Killing form,
so that $\frak k = \frak t \oplus V$.
 The Killing form also induces
$\frak t^* \subset \frak k^*$ from $\frak t \subset \frak k$.
If $F$ is a function on $A$, then by the exponential map, it becomes
a function on $\frak a$. Using the almost complex structure, 
$\frak a^* \cong \frak t^*$. Therefore, the Legendre transform of $F$ can
be written as 
\begin{equation}
L_F : \frak a \longrightarrow \frak t^* .
\label{eq:leg}
\end{equation}

Given $\xi \in \frak k$, we let $\xi^\sharp$ denote its infinitesimal
vector field on $G/N$ induced by the $K$-action.
Let $J$ be the almost complex structure on $G/N$. 
For $\eta = J \xi \in \frak a$, where $\xi \in \frak t$,
we define $\eta^\sharp$ to be the vector field $J \xi^\sharp$.
Let $a \in A \subset KA = G/N$. Then its tangent space is
$T_a(G/N) = \frak h_a^\sharp \oplus V_a^\sharp$.
We recall the
following result from \cite{kn:cg}:

$\spline$

\noindent {\bf Proposition 2.1  } \cite{kn:cg} $ \; $
{\em  Let } $\omega$ {\em be a } $K \times T$ {\em -invariant (1,1)-form
on } $G/N$. {\em Then } $\omega = \sqrt{-1}\partial \bar{\partial}F$,
{\em where } $F \in C^\infty(A)$ {\em by K-invariance. It satisfies }
$\omega(\frak h^\sharp,V^\sharp)_a = 0$.
{\em The K-action is
Hamiltonian, with moment map } $\Phi : G/N \longrightarrow \frak k^*$
{\em satisfying }

{\em (i) } $\;\;$ $\Phi(a) \in \frak t^*$ {\em for all }
$a \in A \subset KA = G/N ;$

{\em (ii) } $\;$ $\Phi : A \longrightarrow \frak t^*$
{\em is given by } $\Phi(a) = \frac{1}{2} L_F(\log a)$.

$\spline$

Let $m = \dim \si, n = \dim \frak t$.
Let $\{ \lambda_1,...,\lambda_r\}$ be the positive roots of $\frak g$,
where $\{ \lambda_1,...,\lambda_n\}$ are simple. 
Here $m \leq n \leq r$. Then $\dim V = 2r$, and
$\dim \frak k = n + 2r$.
In the following proposition, we give a useful decomposition of $V$.
Recall that we define the cell $\si$ in (\ref{eq:cell}) 
using a subset $S$ of the
positive simple roots $\Delta$. 
By switching the roles of $S$ and
$\Delta \backslash S$, we can define another cell $\si'$, with dimension $n-m$. 
We call $\si'$
the complementary cell to $\si$.
Let $J$ be the almost complex
structure on $\frak k \oplus \frak a = \frak g/\frak n$. Recall that,
$V$ is the orthocomplement of $\frak t$ in $\frak k$.

$\spline$

\noindent {\bf Proposition 2.2  }
{\em Let } $\si, \si'$ {\em be complementary cells of dimensions }
$m, n-m$
{\em respectively, where }
$m \leq n \leq r = \frac{1}{2}\dim V$.
{\em There exists a decomposition }
$V = \oplus_1^r V_i$ {\em  into
two dimensional subspaces } $V_i$. {\em  Each } $V_i$
{\em is preserved by J and satisfies } 
$[V_i,V_i] \subset \frak t$. {\em Further, }

{\em (i)} $\;\;\;$ $\frak t_{\si'}^\perp = \oplus_1^m [V_i,V_i] \;\;,$

{\em (ii)} $\;\;$ $\frak t_\si^\perp = \oplus_{m+1}^n [V_i,V_i] \;\; .$

{\em If } $\omega$ {\em is a } $K \times T$ {\em -invariant (1,1)-form on}
$G/N$, {\em then } $\omega(V_i^\sharp,V_j^\sharp)_a = 0$
{\em for all } $i\neq j, a \in A \subset KA = G/N$.

\noindent {\em Proof: }
Let $\{\lambda_1,...,\lambda_r\}$ be the positive roots of $\frak g$,
indexed such that the first $n$ of them are simple.
Further, we can require that
\[ (\lambda_i,\si) > 0 \;,\; (\lambda_i,\si') = 0 \;\;;\;\; i=1,...,m ,\]
and
\[ (\lambda_i,\si) = 0 \;,\; (\lambda_i,\si') > 0 \;\;;\;\; i = m+1,...,n ,\]
where the pairing taken is the Killing form.

Let $\frak g_{\pm i}$ be the root spaces corresponding to $\pm \lambda_i$.
Then there exist $\xi_{\pm i} \in \frak g_{\pm i}$ such that
\begin{equation}
\{ \; \zeta_i = \xi_i - \xi_{-i} \;\;\;,\;\;\;
\gamma_i = \sqrt{-1}(\xi_i + \xi_{-i}) \;\}_{i=1,...,r} 
\label{eq:weyl}
\end{equation}
form a basis of $V$ (\cite{kn:he} p.421). 
Here $\{\zeta_i,\gamma_i\}$ are orthogonal to $\frak t$ because the root
spaces $\frak g_i$ are orthogonal to $\frak h$.
Further, $\{\xi_{\pm i}\}$
can be chosen such that $[\zeta_i,\gamma_i] \in \frak t$, and is dual to 
$\lambda_i \in \frak t^*$ with respect to the Killing form. 
We define 
\[ V_i = {\bf R}(\zeta_i,\gamma_i). \]
Then $[V_i,V_i] \subset \frak t$.
Let $J$ be the almost complex structure on
$\frak k \oplus \frak a = \frak g/\frak n$. 
>From (\ref{eq:weyl}), it follows that
$J$ sends $\zeta_i$ to $\gamma_i$, and sends $\gamma_i$ to $-\zeta_i$.
Therefore, each $V_i$ is preserved by $J$.

For $i=1,...,m$, $(\lambda_i,\si')=0$. Since $[\zeta_i,\gamma_i]$ is dual to
$\lambda_i$, it follows that $[\zeta_i,\gamma_i] \in \frak t_{\si'}^\perp $.
Hence $[V_i,V_i] \subset \frak t_{\si'}^\perp$ for $i=1,...,m$.
But the dual vectors of $\lambda_1,..., \lambda_m$ form a basis of 
$\frak t_{\si'}^\perp$, hence 
$\frak t_{\si'}^\perp = \oplus_1^m [V_i,V_i]$.

For $i=m+1,...,n$, $(\lambda_i,\si)=0$. By similar argument,
$\frak t_\si^\perp = \oplus_{m+1}^n [V_i,V_i]$.

Let $\omega$ be a $K \times T$-invariant {\em (1,1)}-form on $G/N$.
Suppose that $i \neq j$; we want to show that 
$\omega(V_i^\sharp,V_j^\sharp)_a = 0$ for $a \in A \subset KA = G/N$.
Let $p : \frak k \longrightarrow \frak t$ be the 
orthogonal projection, annihilating $V$.
Let $\xi \in V_i, \eta \in V_j$. From (\ref{eq:weyl}), it follows that
$[\xi,\eta]$ is either $0$ or in $V_k$, depending on whether 
$\lambda_i + \lambda_j$ is some positive root $\lambda_k$.
In any case,
\begin{equation}
 p[\xi,\eta]=0 \;\;;\;\; \xi \in V_i,\eta \in V_j .
\label{eq:mes}
\end{equation}
Let $\Phi: G/N \longrightarrow \frak k^*$ be the moment map corresponding
to the $K$-action preserving $\omega$. Then $\Phi(a) \in \frak t^*$,
by Proposition 2.1. Consequently,
\[
\begin{array}{rll}
\omega(\xi^\sharp,\eta^\sharp)_a 
& = (\Phi(a), [\xi,\eta]) & \\
& = (\Phi(a), p [\xi,\eta]) & \mbox{  since  } \Phi(a) \in \frak t^* \\
& = 0 . &
\end{array}
\]
Therefore, $\omega(V_i^\sharp,V_j^\sharp)_a = 0$ for $i \neq j$.
This proves the proposition. \hfill $\Box$

$\spline$

Let $\omega$ be a $K \times T_\si$-invariant Kaehler structure on $G/(P,P)$.
Let $\pi$ be the fibration in (\ref{eq:fib}). 
Then $\pi^* \omega$ is a 
$K \times TA_\si^\perp$-invariant {\em (1,1)}-form
on $G/N$. By Proposition 2.1, it has the form
\[ \pi^* \omega = \sqrt{-1} \partial \bar{\partial} f ,\]
where $f$ is a $K$-invariant function on $G/N$.
Since $G/N = KA$, $f \in C^\infty(A)$.
 We shall show that $f$ can
be replaced with another function $F$ which is in the image of
\[ \pi^* : C^\infty(G/(P,P)) \longrightarrow C^\infty(G/N) , \]
so that we get a potential function for $\omega$.

Let $\si$ be the cell which corresponds to $P$ by (\ref{eq:man}), 
and $\si'$ its complementary cell.
Then $\si'$ defines subgroups $H_{\si'}, T_{\si'}, A_{\si'}$ of $H,T,A$
respectively. 
By taking the orthocomplements of the Lie algebras
$\frak h_{\si'}, \frak t_{\si'}, \frak a_{\si'}$,
we construct the subgroups
$H_{\si'}^\perp, T_{\si'}^\perp, A_{\si'}^\perp$ as before.
Note in particular that
$A = A_\si^\perp A_{\si'}^\perp$.  
Define $F \in C^\infty(A)$ by
\begin{equation}
F = \rho^* f \;\;,\;\;
\rho : A \longrightarrow A_{\si'}^\perp \longrightarrow A \;;
\label{eq:ff}
\end{equation}
where $\rho$ is the composite function of the submersion
$A \longrightarrow A_{\si'}^\perp$ annihilating $A_\si^\perp$, followed by the
inclusion $A_{\si'}^\perp \longrightarrow A$.
By $G/N = KA$, $F$ extends uniquely to be a
$K \times TA_\si^\perp$-invariant function on $G/N$.
Note that $F$ is in the image of $\pi^*$.
We define the $K \times TA_\si^\perp$-invariant {\em (1,1)}-form
\[ \Omega = \sqrt{-1} \partial \bar{\partial}F . \]
We shall show that
\begin{equation}
 \Omega = \pi^* \omega . 
\label{eq:foll}
\end{equation}

Here both $\Omega$ and $\pi^* \omega$ are $K \times TA_\si^\perp$-invariant.
Since $G/N = K A_{\si'}^\perp A_\si^\perp$, we only have to compare them at
$a \in A_{\si'}^\perp$. Also, Proposition 2.1 says that 
$\frak h_a^\sharp$ and $V_a^\sharp$ are complementary with respect to
both $\Omega_a$ and $\pi^* \omega_a$. Therefore, (\ref{eq:foll}) will follow
if we can show that
\begin{equation}
 \Omega(\xi^\sharp,\eta^\sharp)_a 
= \pi^* \omega(\xi^\sharp,\eta^\sharp)_a \;\;;\;\; 
\xi,\eta \in \frak h \mbox{  or  } \xi,\eta \in V \,,\, a \in A_{\si'}^\perp. 
\label{eq:cruc}
\end{equation}
This will be checked by the following two lemmas. 
Recall that $L_F,L_f : \frak a \longrightarrow \frak t^*$ are the Legendre
transforms of $F$ and $f$, described in (\ref{eq:leg}).

$\spline$

\noindent {\bf Lemma 2.3  } 
$ \Omega(\xi^\sharp,\eta^\sharp)_a = \pi^* \omega(\xi^\sharp,\eta^\sharp)_a $
{\em  for all } $ \xi,\eta \in V , a \in A_{\si'}^\perp$.

\noindent {\em Proof: }
By Proposition 2.2, the spaces 
$(V_1)_a^\sharp,...,(V_r)_a^\sharp$ are pairwise complementary with respect
to $\Omega_a$ and $\pi^* \omega_a$, $a \in A_{\si'}^\perp$. Therefore,
to prove the statement in this lemma, we may consider 
$\xi, \eta \in V_i$ for each component $V_i$ seperately. 
Since each $V_i$ is two dimensional, it suffices to consider 
$\xi = \zeta_i, \eta = \gamma_i$. Let 
\[ \Phi_F, \Phi_f : G/N \longrightarrow \frak k^* \]
be the moment maps of the $K$-actions preserving $\Omega, \pi^* \omega$
respectively. We recall from Proposition 2.1 that
$\Phi_F(a)=\frac{1}{2}L_F(\log a), \Phi_f(a)=\frac{1}{2}L_f(\log a)$.
We follow the indices $i=1,...,r$ used in Proposition 2.2, as well as
the cells $\si, \si'$ of dimensions $m, n-m$ respectively. In what follows, we
break up our arguments into three cases, according to the different values
of the index $i$.

\noindent {\em Case 1: } $i=1,...,m$.
\[
\begin{array}{rl}
\Omega(\zeta_i^\sharp,\gamma_i^\sharp)_a 
& = (\Phi_F(a),[\zeta_i,\gamma_i]) \\
& = (\frac{1}{2} L_F (\log a), [\zeta_i,\gamma_i]) .
\end{array}
\]
By Proposition 2.2, $[\zeta_i,\gamma_i] \in \frak t_{\si'}^\perp$, 
for $i=1,...,m$.
By (\ref{eq:ff}), $L_F(\log a)$ and $L_f(\log a)$ 
agree on $\frak t_{\si'}^\perp$, for $a \in A_{\si'}^\perp$.
Therefore, the last expression is
\[ 
\begin{array}{rl}
(\frac{1}{2} L_f(\log a), [\zeta_i,\gamma_i]) 
& = (\Phi_f(a), [\zeta_i,\gamma_i]) \\
& = \pi^*\omega(\zeta_i^\sharp,\gamma_i^\sharp)_a .
\end{array}
\]

\noindent {\em Case 2: } $i=m+1,...,n$.

We recall (\ref{eq:weyl}), which implies that
\begin{equation}
[v,\zeta_i] = \sqrt{-1}(\lambda_i,v) \gamma_i \;\;,\;\;
[v,\gamma_i] = -\sqrt{-1}(\lambda_i,v) \zeta_i
\label{eq:ei}
\end{equation}
for all $v \in \frak t$. Therefore, the Lie algebra
$\frak k^\si$ of $K^\si$ is given by
\[ \frak k^\si \; = \; \{\xi \in \frak k\;;\; [\xi,\si] = 0\} \; = \;
\frak t \oplus_{(\lambda_i,\si)=0}V_i . \]
The center of this Lie algebra is $\frak t_\si$, hence the semi-simple
Lie algebra $\frak k_{ss}^\si$ is given by
\begin{equation}
\frak k_{ss}^\si \; = \; \frak t_\si^\perp \oplus_{(\lambda_i,\si)=0} V_i .
\label{eq:kcss}
\end{equation}
For $i=m+1,...,n$, $(\lambda_i,\si)=0$; hence
$\zeta_i,\gamma_i \in \frak k_{ss}^\si$. But $K_{ss}^\si$
is in the fiber of $\pi$, so $\imath(\xi^\sharp)\pi^* \omega_a = 0$
for all $\xi \in V_i$. 

We shall show that
\[ \imath(\xi^\sharp) \Omega_a = 0 \]
for all $\xi \in V_i$.
Since each $V_i$
is two dimensional, this will follow if we can show that
$\Omega(\zeta_i^\sharp,\gamma_i^\sharp)_a = 0$, for $i=m+1,...,n$. But
\[ \Omega(\zeta_i^\sharp,\gamma_i^\sharp)_a 
= (\frac{1}{2} L_F(\log a), [\zeta_i,\gamma_i]) =0 ,\]
since $[\zeta_i,\gamma_i] \in \frak t_\si^\perp$ and 
by (\ref{eq:ff}), $L_F(\log a)$ vanishes there.

\noindent {\em Case 3: } $i=n+1,...,r$.

>From Cases 1, 2, we see that $L_F(\log a), L_f(\log a) \in \frak t^*$
agree on the spaces $\frak t_\si^\perp, \frak t_{\si'}^\perp$. Since
$\frak t = \frak t_\si^\perp \oplus \frak t_{\si'}^\perp$,
it follows that $L_F(\log a) = L_f(\log a) \in \frak t^*$. Therefore,
\[
\begin{array}{rl}
\Omega(\zeta_i^\sharp,\gamma_i^\sharp)_a
& = (\Phi_F(a), [\zeta_i,\gamma_i]) \\
& = (\frac{1}{2} L_F(\log a), [\zeta_i,\gamma_i]) \\
& = (\frac{1}{2} L_f(\log a), [\zeta_i,\gamma_i]) \\
& = (\Phi_f(a), [\zeta_i,\gamma_i]) \\
& = \pi^* \omega(\zeta_i^\sharp,\gamma_i^\sharp)_a .
\end{array}
\]

This proves Lemma 2.3. \hfill $\Box$

$\spline$

\noindent {\bf Lemma 2.4  }
$ \Omega(\xi^\sharp,\eta^\sharp)_a 
= \pi^* \omega (\xi^\sharp,\eta^\sharp)_a $
{\em  for all } $\xi,\eta \in \frak h, a \in A_{\si'}^\perp .$

\noindent {\em Proof: } Let $\frak h_\si, \frak h_{\si'}$ 
denote the subalgebras
of $\frak h$, by taking the complex linear spans of $\si, \si'$ respectively.
Let $\frak h_\si^\perp, \frak h_{\si'}^\perp$ denote their orthocomplements
with respect to the Killing form.
Then $\frak h = \frak h_\si^\perp \oplus \frak h_{\si'}^\perp$.

\noindent {\em Case 1: } $\xi,\eta \in \frak h_{\si'}^\perp$.

Let $\iota : H_{\si'}^\perp \longrightarrow H$ denote the inclusion.
>From (\ref{eq:ff}), we get
\[ \sqrt{-1} \partial \bar{\partial} (\iota^* F)
= \sqrt{-1}\partial \bar{\partial} (\iota^* f) , \]
where $\partial, \bar{\partial}$ are Dolbeault operators on $H_{\si'}^\perp$
here. Therefore, given $a \in A_{\si'}^\perp \subset H_{\si'}^\perp$,
\[
\begin{array}{rl}
\Omega(\xi^\sharp,\eta^\sharp)_a
& = (\iota^* \Omega)(\xi^\sharp,\eta^\sharp)_a \\
& = (\sqrt{-1}\partial \bar{\partial}(\iota^* F))(\xi^\sharp,\eta^\sharp)_a \\
& = (\sqrt{-1}\partial \bar{\partial}(\iota^* f))(\xi^\sharp,\eta^\sharp)_a \\
& = (\iota^* \pi^* \omega)(\xi^\sharp,\eta^\sharp)_a \\
& = \pi^* \omega(\xi^\sharp,\eta^\sharp)_a .
\end{array}
\]

\noindent {\em Case 2: } $\xi \in \frak h_\si^\perp$.

We shall show that
\begin{equation}
 \imath(\xi^\sharp) \pi^* \omega_a 
\;=\; \imath(\xi^\sharp) \Omega_a \;=\; 0 ,
\label{eq:tir}
\end{equation}
which completes the proof of this lemma. 
Since $\pi^* \omega$ and $\Omega$ are {\em (1,1)}-forms,
it suffices to check (\ref{eq:tir}) for $\xi \in \frak t_\si^\perp$.

The fiber of $\pi$ is $K_{ss}^\si \times A_\si^\perp$, which contains
$H_\si^\perp$. Therefore, 
\[ \imath(\xi^\sharp) \pi^* \omega_a = 0 . \]

We observe that, as complex manifolds,
\[ H = {\bf C}^n/{\bf Z}^n \;\;,
H_\si^\perp = {\bf C}^{n-m}/{\bf Z}^{n-m} \;\;,
H_{\si'}^\perp = {\bf C}^m/{\bf Z}^m , \]
and $H = H_\si^\perp H_{\si'}^\perp$.
We introduce complex coordinates $\{z_1,...,z_m\}$ on $H_{\si'}^\perp$
as well as
$\{z_{m+1},...,z_n\}$ on $H_\si^\perp$; so that $H$ adopts the product
coordinates. Let $z = x + \sqrt{-1}y$, and we let $x, y$ be the coordinates
on $T, A$ respectively. 
>From $H=TA,G/N=KA$ and $T \subset K$, we get a natural holomorphic imbedding
$\iota : H \longrightarrow G/N$.
Then $\iota^* F$, being $T$-invariant, 
is a function on $y$ only. For simplicity we still denote it as $F$.
It follows from (\ref{eq:ff}) that
\[ \frac{\partial F}{\partial y_i} = 0 
\mbox{  for  } i=m+1,...,n.\]
Therefore, for $a \in A_{\si'}^\perp$,
\begin{equation}
\begin{array}{rl}
\imath(\xi^\sharp)(\iota^* \Omega)_a
& = \imath(\xi^\sharp)(\sqrt{-1}\partial \bar{\partial}F)_a \\
& = \imath(\xi^\sharp)(\frac{1}{2} \sum_{j,k=1}^n
\frac{\partial^2 F}{\partial y_j \partial y_k} dx_j \wedge dy_k) \\
& = \imath(\xi^\sharp)(\frac{1}{2} \sum_{j,k=1}^m
\frac{\partial^2 F}{\partial y_j \partial y_k} dx_j \wedge dy_k) .
\end{array}
\label{eq:qwe}
\end{equation}
On the other hand, since $\xi \in \frak t_\si^\perp$, 
the vector field $\xi^\sharp$ on $H$
is of the form
\[ \xi^\sharp = \sum_{m+1}^n c_i \frac{\partial}{\partial x_i} . \]
This, together with (\ref{eq:qwe}), imply that
\[ \imath(\xi^\sharp) \Omega_a = 0 .\]
This proves (\ref{eq:tir}).
Combining the results in Cases 1,2, we have proved Lemma 2.4. \hfill $\Box$

$\spline$

Lemmas 2.3 and 2.4 imply (\ref{eq:cruc}),
and hence (\ref{eq:foll}). Namely, we have shown that 
given a $K \times T_\si$-invariant
Kaehler structure $\omega$ on $G/(P,P)$, there exists a function $F$,
which is in the image of $\pi^*$ by virtue of (\ref{eq:ff}), such that
\[ \pi^* \omega = \sqrt{-1}\partial \bar{\partial} F .\]
Since $F$ is in the image of $\pi^*$,
 and since $\pi^*$ is injective, it follows that
$\omega$ itself has a potential function.

Conversely, suppose that a $K$-invariant Kaehler structure $\omega$
on $G/(P,P)$ has a potential function $F$. Averaging by the compact group $K$
if necessary, we may assume that $F$ is $K$-invariant. But by (\ref{eq:nnsc}),
this means that $F$ is just a function on $A_\si$, and is automatically
$K \times T_\si$-invariant. Then $\omega$ is also $K \times T_\si$-invariant.
This proves Theorem I.

We note that our arguments do not require $\omega$ to be 
positive definite. Namely,
Theorem I holds even if $\omega$ is merely a $K$-invariant
{\em (1,1)}-form. In the next section, we use the moment map to derive
a necessary and sufficient condition for a $K \times T_\si$-invariant
{\em (1,1)}-form to be Kaehler.

\newpage
\begin{center}
\section{MOMENT MAP}
\end{center}
\setcounter{equation}{0}

Let $\omega$ be a $K \times T_\si$-invariant {\em (1,1)}-form on $G/(P,P)$, 
with moment map 
\[ \Phi : G/(P,P) \longrightarrow \frak k^* \]
corresponding to the Hamiltonian
action of $K$ on $G/(P,P)$ preserving $\omega$.
It is easy to see that this action is Hamiltonian; either from the
semi-simplicity of $K$ (\cite{kn:gs2}, \S 26), or from the fact that
$\omega = \sqrt{-1} \partial \bar{\partial} F$ implies $\omega = d \beta$
for some $K$-invariant real 1-form $\beta$ (\cite{kn:am}, Theorem 4.2.10).
We shall study the moment map $\Phi$, and derive a necessary and sufficient 
condition for $\omega$ to be Kaehler.

Suppose now that $\omega$ is a $K \times T_\si$-invariant Kaehler structure.
We want to derive the two conditions stated in Theorem II.
By Theorem I, $\omega$ has a potential function $F$.
Averaging by $K$ if necessary, we may assume that $F$ is $K$-invariant.
By (\ref{eq:nnsc}),
$G/(P,P)= (K/K_{ss}^\si)\times A_\si$;
so the $K$-invariant function $F$ is just a function on $A_\si$.
Let $\pi$ be the fibration in (\ref{eq:fib}). Then
\[ \Phi \circ \pi : G/N \longrightarrow \frak k^* \]
is the moment map corresponding to the $K$-action on
$(G/N, \pi^* \omega)$. Recall that $P$ corresponds to a cell $\si$
via (\ref{eq:man}).
Also, $G/N = KA$ and 
$G/(P,P) = (K/K_{ss}^\si) \times A_\si$ induce the inclusions
\[ A \hookrightarrow \{e\}\times A \subset KA = G/N \;\;,\;\;
A_\si \hookrightarrow \{e K_{ss}^\si\} \times A_\si 
\subset (K/K_{ss}^\si)\times A_\si =G/(P,P) \;.\]
Therefore, we can regard $A$
and $A_\si$ as contained in $G/N$ and $G/(P,P)$ respectively.
Note that $\pi(A)=A_\si$.
>From Proposition 2.1, we see that
\[ (\Phi \circ \pi)(A) \subset \frak t^* .\]
Since the fibration $\pi$ sends $A$ to $A_\si$, it follows that
$\Phi(A_\si) \subset \frak t^*$. 
By $K$-equivariance of $\Phi$, $\Phi |_{A_\si}$ determines $\Phi$
entirely. The exponential map from $\frak a_\si$ to $A_\si$ is a
diffeomorphism, and we let $\log$ be its inverse. 
This way, the potential function $F$ becomes a function on $\frak a_\si$.
Then, by the almost complex structure, $\frak a_\si^* \cong \frak t_\si^*$.
Consequently, the Legendre transform of $F$ is
\[ L_F : \frak a_\si \longrightarrow \frak t_\si^* . \]
We shall show that
\[ \Phi : A_\si \longrightarrow \frak t^* \]
is given by $\Phi(a) = \frac{1}{2} L_F(\log a)$ for all $a \in A_\si$. Let 
\[ \imath : H_\si \longrightarrow G/(P,P) \]
be the natural holomorphic imbedding of $H_\si= T_\si A_\si$. 
Then $\imath^* \omega$ is a $T_\si$-invariant
Kaehler structure on $T_\si A_\si$, with potential function $\imath^* F$.
For simplicity, we still write $\imath^* F$ as $F$.
Let $m$ be the dimension of the cell $\si$.
Then, as a complex manifold, $H_\si = {\bf C}^m/{\bf Z}^m$. Therefore,
we can introduce complex coordinates $\{z_1,...,z_m\}$ on $H_\si$, where 
\begin{equation}
\begin{array}{c}
H_\si = {\bf C}^m/{\bf Z}^m = \{z_1,...,z_m\} \;,\;
T_\si = {\bf R}^m/{\bf Z}^m = \{x_1,...,x_m\} \;,\\
A_\si = {\bf R}^m = \{y_1,...,y_m\} \;,\; z_i = x_i + \sqrt{-1} y_i \;.
\label{eq:cz}
\end{array}
\end{equation}
Since $F$ is $T_\si$-invariant, it is a function on $y$ only.
Then $\imath^*\omega$ becomes (here $\partial,\bar{\partial}$ are
Dolbeault operators on $H_\si$)
\begin{equation}
\imath^* \omega = \sqrt{-1} \partial \bar{\partial} F
= \frac{1}{2} \sum_{j,k=1}^m 
\frac{\partial^2 F}{\partial y_j \partial y_k} dx_j \wedge dy_k ,
\label{eq:small}
\end{equation}
where $F \in C^\infty({\bf R}^m)$. Since $\omega$ is Kaehler, 
so is $\imath^* \omega$; and (\ref{eq:small}) says that $\imath^*\omega$
is Kaehler if and only if the Hessian matrix of $F$ is positive definite,
i.e. $F$ is strictly convex.

The moment map $\Phi$ of the $K$-action on
$(G/(P,P),\omega)$ restricts to be the moment map 
$\Phi '$ of the $T_\si$-action on
$(T_\si A_\si, \imath^*\omega)$. Let
\[ \beta = -\frac{1}{2} \sum_{j=1}^m
\frac{\partial F}{\partial y_j} dx_j \]
be a $T_\si$-invariant 1-form on $T_\si A_\si$. From (\ref{eq:small}),
it follows that $d\beta = \imath^* \omega$, so the moment map $\Phi '$ of 
the $T_\si$-action is
\[
\begin{array}{rl}
(\Phi ' (ta) , \xi) & = - (\beta, \xi^\sharp)(ta) \\
& = (\frac{1}{2} \sum_{j=1}^m \frac{\partial F}{\partial y_j} dx_j, 
\sum_{k=1}^m \xi_k \frac{\partial}{\partial x_k})(ta) \\
& = \frac{1}{2} \sum_{j=1}^m \frac{\partial F}{\partial y_j}(a) \xi_j \\
& = \frac{1}{2} (L_F(a), \xi) \;,
\end{array}
\]
where $ta \in T_\si A_\si, \xi \in \frak t = {\bf R}^m$. 
Our computation identifies $\frak a$ with $A$ by the exponential map,
so in fact $\Phi '(ta) = \frac{1}{2} L_F (\log a)$ 
for all $ta \in T_\si A_\si$.
But $\Phi$ and $\Phi '$ agree on $A_\si$, so 
$\Phi(a) = \frac{1}{2} L_F(\log a)$. 
Hence $\Phi(A_\si) \subset \frak t_\si^*$.
We claim further that $\Phi(A_\si) \subset \si$:

Let $V_i \subset V \subset \frak k$ be the subspaces constructed in 
Proposition 2.2, and let $\{\zeta_i, \gamma_i\} \in V_i$ be
the vectors in (\ref{eq:weyl}). 
Recall that these indices are made with respect to the positive roots
$\{\lambda_i\}$. 
Since $G/(P,P) = (K/K_{ss}^\si)\times A_\si$, the infinitesimal
vector fields $\zeta_i^\sharp,\gamma_i^\sharp$ 
on $G/(P,P)$ are non-zero if and only if
$\zeta_i,\gamma_i \not{\!\in} \frak k_{ss}^\si$. 
By (\ref{eq:kcss}), this is
equivalent to $(\lambda_i,\si) > 0$.
Let $J$ be the almost
complex structure in $G/(P,P)$, $a \in A_\si$, and $(\lambda_i,\si) > 0$
so that $\zeta_i^\sharp, \gamma_i^\sharp \neq 0$. By (\ref{eq:weyl}), 
$J \zeta_i = \gamma_i$.
Since $\omega$ is Kaehler,
\begin{equation}
\begin{array}{rl}
0 & < \omega(\zeta_i^\sharp, J \zeta_i^\sharp)_a \\
& = \omega(\zeta_i^\sharp, \gamma_i^\sharp)_a \\
& = (\Phi(a), [\zeta_i,\gamma_i]) \\
& = (\Phi(a), \lambda_i) .
\end{array}
\label{eq:moc}
\end{equation}
We have shown that, for all $a \in A_\si$, $(\Phi(a),\lambda_i) > 0$
whenever $\lambda_i$ is a positive root satisfying $(\lambda_i,\si) > 0$.
This, together with $\Phi(A_\si) \subset \frak t_\si^*$, imply that
$\Phi(A_\si) \subset \si$, as claimed.

We have shown that if $\omega$ is Kaehler, then the two conditions stated
in Theorem II have to be satisfied. We next show that, conversely, 
these two conditions are sufficient for $\omega$ to be Kaehler.

Recall that
the infinitesimal vector field $\xi^\sharp$ on $G/(P,P)$ vanishes if 
$\xi \in \frak k_{ss}^\si$. Hence the tangent space at 
$a \in A_\si \subset G/(P,P)$ is spanned by 
$(\frak k_{ss}^{\si \perp})_a^\sharp, (\frak a_\si)_a^\sharp$.
Here we define $\eta^\sharp$ for $\eta = J \xi \in \frak a_\si$
by $\eta^\sharp = J \xi^\sharp$, where $\xi \in \frak t_\si$.
However, it follows from (\ref{eq:kcss}) that
\[\frak k_{ss}^{\si \perp} \; = \; 
\frak t_\si \oplus_{(\lambda_i,\si)>0} V_i , \]
where $V_i$ is the space described in Proposition 2.2.
Here the distinct $V_i$ are orthogonal to one another, due to the
orthogonality of the root spaces $\frak g_i$ (\cite{kn:he} p.166).
Consequently, the tangent space at $a \in A_\si \subset G/(P,P)$ is
\begin{equation}
T_a(G/(P,P))= (\frak h_\si)_a^\sharp \oplus_{(\lambda_i,\si)>0}(V_i)_a^\sharp .
\label{eq:sep}
\end{equation}

We claim that 
$\omega(\frak h_\si^\sharp, V_i^\sharp)_a =
\omega(V_i^\sharp, V_j^\sharp)_a = 0$, for $i \neq j$:

Since $J$ preserves $\frak h_\si$ and $V_i$, and $\omega$ is a 
{\em (1,1)}-form,
the first part follows if we can
show that $\omega(\frak t_\si^\sharp, V_i^\sharp)_a = 0$.
Let $p : \frak k \longrightarrow \frak t$ be the orthogonal projection,
annihilating $V$. Let $\xi \in \frak t_\si, \eta \in V_i$. Then 
$p[\xi,\eta]= 0$, by (\ref{eq:ei}).
Since $\Phi(a) \in \frak t^*$ for $a \in A$,
\[\omega(\xi^\sharp,\eta^\sharp)_a \;=\; (\Phi(a),[\xi,\eta])
\;=\; (\Phi(a),p[\xi,\eta]) \;=\; 0 .\]
Hence $\omega(\frak h_\si^\sharp,V_i^\sharp)_a = 0$.
For $i \neq j$, it follows from (\ref{eq:mes}) that
$p[V_i,V_j]=0$. So, by similar argument, 
$\omega(V_i^\sharp, V_j^\sharp)_a = 0$ as claimed.
 
Therefore, by $K$-invariance of $\omega$ and (\ref{eq:sep}),
the positive definite condition of $\omega$
follows if we can check that
\begin{equation}
\omega(\xi^\sharp,J \xi^\sharp)_a > 0 \;\;;\;\;
\xi \in \frak h_\si \mbox{  or  } \xi \in V_i \,,\,
(\lambda_i,\si)>0, a \in A_\si .
\label{eq:pd}
\end{equation}

But they follow from the two conditions of Theorem II:
Condition (i) of Theorem II
implies that the expression in (\ref{eq:small}) is positive definite and hence
(\ref{eq:pd}) holds for $\xi \in \frak h_\si$.
Condition (ii) of Theorem II
implies that $(\Phi(a),\lambda_i)>0$ whenever $(\lambda_i,\si)>0$,
so it follows from (\ref{eq:moc}) that
(\ref{eq:pd}) holds
for $\xi \in V_i$.
This proves Theorem II.

\newpage
\begin{center}
\section{LINE BUNDLE}
\end{center}
\setcounter{equation}{0}

Fix a $K \times T_\si$-invariant Kaehler structure $\omega$ on $G/(P,P)$.
By Theorem I, $\omega$ has a potential function
$F$. Recall that $P$ determines the subgroup $A_\si$ by (\ref{eq:man}).
By $K$-invariance and (\ref{eq:nnsc}),
 we can regard $F$ as a function on $A_\si$.
In particular, the expression 
$\omega = \sqrt{-1}\partial \bar{\partial}F$ also implies that $\omega$
is exact. Hence $\omega$ is integral, and there exists a complex line bundle
${\bf L}$ on $G/(P,P)$
whose Chern class is $[\omega]=0$, equipped with a connection
$\nabla$ whose curvature is $\omega$,
as well as an invariant Hermitian structure $<,>$ (\cite{kn:gs}, \cite{kn:ko}).
The line bundle
${\bf L}$ is trivial since $[\omega]=0$, but the connection $\nabla$
gives rise to interesting geometry. We say that a section $s$ is holomorphic
if $\, \nabla s \,$ annihilates anti-holomorphic vector fields on $G/(P,P)$.
We shall show that the $K \times T_\si$-action on $G/(P,P)$ lifts to a
$K \times T_\si$-representation on the space of holomorphic sections 
of ${\bf L}$. 
To do this, we shall construct a global trivialization of ${\bf L}$. 
The following topological property of $G/(P,P)$ is useful in
this construction:
 
$\spline$

\noindent {\bf Lemma 4.1} \hspace{2mm} $H^{1}(G/(P,P), {\bf C}) = 0$.

\noindent {\em Proof:  }  
By (\ref{eq:nnsc}), $G/(P,P) = (K/K_{ss}^\si) \times A_\si$.
Since $A_\si$ is Euclidean, it suffices to show that
 $H^{1} (K/K_{ss}^\si, {\bf C}) = 0$.

The fibration $K \longrightarrow K/K_{ss}^\si$ 
 induces a long exact sequence of homotopy groups,
\begin{equation}
... \longrightarrow {\pi}_{1}(K) \longrightarrow
{\pi}_{1}(K / K_{ss}^\si) \longrightarrow {\pi}_{0}(K_{ss}^\si)
\longrightarrow ...  \label{eq:longex}
\end{equation}
However, by (\cite{kn:bd} p.223),
\[ {\pi}_{1}(K) \; \cong \;
\mbox{ker}(\exp : \frak t \rightarrow T)/{\bf Z}(\mbox{roots of }\frak k) .\]
Therefore, since $K$ is semi-simple, ${\pi}_{1}(K)$ is finite.
 By compactness of $K_{ss}^\si$, ${\pi}_{0}(K_{ss}^\si)$ is finite.
Hence ${\pi}_{1}(K/K_{ss}^\si)$, being caught in the middle
in (\ref{eq:longex}), is also finite.
It follows that 
\[ H^{1}(K/K_{ss}^\si, {\bf C}) \cong
 Hom({\pi}_{1}(K/K_{ss}^\si), {\bf C}) = 0, \]
which proves the lemma. \hfill $\Box$

$\spline$

We return to our pre-quantum line bundle ${\bf L}$ on $G/(P,P)$,
corresponding to the $K \times T_\si$-invariant Kaehler structure $\omega$.
Let $\beta$ be the 1-form $- \sqrt{-1} \partial F$, so 
$ d \beta = \omega . $
We claim that

$\spline$

\noindent {\bf Proposition 4.2} {\em There exists a non-vanishing section}
$s_o$ {\em on } ${\bf L}$ {\em  , with the property}
\begin{equation}
\beta = \frac{1}{\sqrt{-1}} \frac{\nabla s_o}{s_o}. \label{eq:copi}
\end{equation}
{\em This section is 
unique up to a non-zero constant multiple, and is holomorphic.
Up to a non-zero constant,}
\[ <s_o,s_o> = e^{-F} . \]

\noindent {\em Proof: } Since $[\omega] = 0$, ${\bf L}$ is a trivial bundle;
so there exists a nowhere zero section $s_1$ of ${\bf L}$. Let
\[ \alpha = \frac{1}{\sqrt{-1}} \frac{\nabla s_1}{s_1} . \]
By the definition of the curvature form on ${\bf L}$, $d \alpha = \omega$;
so $d (\beta - \alpha) = 0$. Since 
$H^{1}(G/(P,P), {\bf C}) = 0$, there exists a 
complex-valued function $f$ such that $\beta = \alpha + df$. 
Let $s_o = (\exp \sqrt{-1} f) s_1$. Then
\[ \frac{1}{\sqrt{-1}} \frac{\nabla s_o}{s_o} = 
\frac{1}{\sqrt{-1}} \frac{\nabla s_1}{s_1} + df = \beta .\]
This proves the existence of a holomorphic section $s_o$ satisfying
(\ref{eq:copi}).

Suppose that
$s_{1}$ and $s_{2}$ are two sections satisfying this formula. Let
$h = \frac{s_{2}}{s_{1}} .$ Then
\[ \frac{1}{\sqrt{-1}} \frac{\nabla s_{2}}{s_{2}} =
\frac{1}{\sqrt{-1}} \frac{\nabla s_{1}}{s_{1}} + 
\frac{1}{\sqrt{-1}} d \log h ,\]
which implies that $h$ is a constant. Hence, up to a constant, the solution
of (\ref{eq:copi}) is unique.

If $v$ is an anti-holomorphic vector field, then
\[ \frac{1}{\sqrt{-1}} \frac{{\nabla}_{v} s_o}{s_o} = \iota (v) \beta = 0 , \]
as $\beta$ is a form of type {\em (1,0)}.
 Hence $s_o$ is holomorphic. 
Since $\beta$ is $K \times T_\si$-invariant,
$s_o$ induces a $K \times T_\si$-representation on the space of holomorphic
sections on ${\bf L}$, where $s_o$ is $K \times T_\si$-invariant.
Namely, given a
holomorphic section $f s_o$ of ${\bf L}$ (note $s_o$ is non-vanishing),
$K \times T_\si$ acts by
\begin{equation}
 L_k^* R_t^*(f s_o) \, = \, (L_k^* R_t^* f)s_o \;\;;\;\; 
k \in K, t \in T_\si ,
\label{eq:kt}
\end{equation}
where $L_k^* R_t^* f$
denotes the standard action on the holomorphic functions lifted from 
the $K \times T_\si$-action on $G/(P,P)$. Hence
 $s_o$ defines a $K \times T_\si$-equivariant trivialization.

For this section $s_o$, we now show that $<s_o,s_o> = e^{-F}$. By
$K$-invariance, it suffices to show that this is the case when restricted to
$A_\si$. Let $\si$ be the cell corresponding to the parabolic subgroup $P$,
and let $m$ be the dimension of $\si$.
We write
$ T_\si A_\si = {\bf C}^{m} / {\bf Z}^{m} = \{ z_1,...,z_m\} $
as in (\ref{eq:cz}), so that $F$, being $T_\si$-invariant,
is a function on $y$ only.
Let $\imath : T_\si A_\si \longrightarrow G/(P,P)$ be the 
natural inclusion. Then
\begin{equation}
{\imath}^{*} \beta = - \sqrt{-1} \partial F
= \frac{1}{2} \sum_1^m \frac{\partial F}{\partial y_{i}} d z_{i}.
\label{eq:tigalah}
\end{equation}
Let ${\nabla}_{i} = {\nabla}_{\frac{\partial}{\partial y_{i}}}$. Then 
\[ \frac{\partial}{\partial y_{i}} <s_o,s_o> 
= \hspace{2mm} <{\nabla}_{i} s_o,s_o> + <s_o, {\nabla}_{i} s_o>. \]
However, by (\ref{eq:copi}) and (\ref{eq:tigalah}), 
\[ \frac{{\nabla}_{i} s_o}{s_o} = 
\sqrt{-1} ( \beta, \frac{\partial}{\partial y_{i}}) =
- \frac{1}{2} \frac{\partial F}{\partial y_{i}} \]
so 
\[ \frac{\partial}{\partial y_{i}} \log <s_o, s_o> =
- \frac{\partial F}{\partial y_{i}}. \]
Therefore, up to a non-zero constant multiple,
\[ <s_o, s_o> = e^{-F} . \]
This proves the proposition. \hfill $\Box$

$\spline$

Let ${\cal O}({\bf L})$ denote the space of holomorphic sections of 
the line bundle ${\bf L}$ on $G/(P,P)$.
By Proposition 4.2, $s_o$ induces a $K \times T_\si$-representation on
${\cal O}({\bf L})$,
given by (\ref{eq:kt}), where $s_o$ is $K \times T_\si$-invariant.
 Let $\lambda \in \frak t_\si^*$ be a dominant integral weight, and
let ${\cal O}({\bf L})_\lambda$ denote the holomorphic sections that
transform by $\lambda$ under the right $T_\si$-action.
Since this action commutes with the left $K$ action, 
${\cal O}({\bf L})_\lambda$ is a $K$-subrepresentation of ${\cal O}({\bf L})$.
We now show that the $K$-finite vectors in ${\cal O}({\bf L})$
decompose into $\{{\cal O}({\bf L})_\lambda \;;\; \lambda \in \overline{\si}\}$ 
as irreducible $K$-representations with highest weights $\lambda$. 
Using the holomorphic section $s_o$
of Proposition 4.2, it suffices to consider the holomorphic functions
${\cal O}(G/(P,P))$; since
\[ {\cal O}(G/(P,P)) \otimes s_o = {\cal O}({\bf L}) \]
is a $K \times T_\si$-equivariant trivialization.

Recall that $\overline{W}$ is the closure of the Weyl chamber $W$, and 
$\si \subset \overline{W}$ is the cell corresponding to $P$. 
Let $\overline{\si}$ denote its closure in $\overline{W}$.
For a dominant integral weight $\lambda \in \frak t^*$,
let ${\cal O}_\lambda \subset {\cal O}(G/(P,P))$ denote the holomorphic
functions that transform by $\lambda$ under the right $T_\si$-action.
Since the right $T_\si$-action commutes with the left $K$-action,
each ${\cal O}_\lambda$ is a $K$-representation space.

$\spline$

\noindent {\bf Proposition 4.3  }  
{\em  The irreducible K-representation with highest weight } $\lambda$
{\em occurs in } ${\cal O}(G/(P,P))$
{\em if and only if } $\lambda \in \overline{\si}$.
{\em For } $\lambda \in \overline{\si}$, {\em it occurs with multiplicity one,
and is given by } ${\cal O}_\lambda$.

\noindent {\em Proof: }
The fibration $\pi$ of (\ref{eq:fib})
induces an injection of holomorphic functions,
\[ \pi^* : {\cal O}(G/(P,P)) \longrightarrow {\cal O}(G/N) .\]
This map intertwines with the $K \times T_\si$-action.

Let $\lambda$ be a dominant integral weight, but suppose that
$\lambda \not{\!\!\in} \overline{\si}$. We shall show that the $K$-irreducible
with highest weight $\lambda$ does not occur in ${\cal O}(G/(P,P))$.
By the Borel-Weil theorem, the $K$-irreducible with highest weight $\lambda$
occurs in ${\cal O}(G/N)$ with multiplicity one, and can be taken as the
holomorphic functions in $G/N$ that transform by $\lambda$ under the right 
$T$-action. We denote this space by 
$V_\lambda \subset {\cal O}(G/N)$. Since
$\pi^*$ is injective, it suffices
to show that
\begin{equation}
\pi^* {\cal O}(G/(P,P)) \cap V_\lambda = 0 .
\label{eq:koso}
\end{equation}
Since $\lambda \not{\!\in} \overline{\si}$, $(\lambda, \xi) \neq 0$ for some
$\xi \in \frak t_\si^\perp$. Let $0 \neq f \in V_\lambda$. Then 
the right action $R_\xi^*$ on $V_\lambda$ satisfies
\begin{equation}
R_\xi^* f = (\lambda, \xi) f \neq 0 .
\label{eq:cmd}
\end{equation}
Since $T_\si^\perp$ is in the fiber of $\pi$, the image of $\pi^*$ is 
$T_\si^\perp$-invariant.
Therefore, (\ref{eq:cmd}) says that $f$ cannot be in the image of $\pi^*$.
This proves (\ref{eq:koso}). 

Conversely, suppose that
$\lambda \in \overline{\si}$ is a dominant integral weight.
We again let $V_\lambda \subset {\cal O}(G/N)$ be the holomorphic functions
that transform by $\lambda$. By the Borel-Weil theorem, $V_\lambda$ is
an irreducible representation with highest weight $\lambda$,
and such irreducible occurs with multiplicity one.
 Therefore, 
to complete the proof of Proposition 4.3,
we need to show
\begin{equation}
V_\lambda \subset \pi^*{\cal O}(G/(P,P)) \;,\; \lambda \in \overline{\si}.
\label{eq:tt}
\end{equation}

Recall from (\ref{eq:reyer}) that the fiber of $\pi$ is 
$(K_{ss}^\si)_{\bf C}/(M \cap N) = K_{ss}^\si \times A_\si^\perp$.
Choose a fiber of $\pi$, and let
\[ \imath : K_{ss}^\si A_\si^\perp \hookrightarrow G/N \]
be a holomorphic $K_{ss}^\si \times A_\si^\perp$-equivariant 
imbedding as this fiber of $\pi$.
Let $f \in V_\lambda$. We claim that $f$ is constant on this fiber:

By applying the Borel-Weil theorem on
$(K_{ss}^\si)_{\bf C}/(M \cap N) = K_{ss}^\si \times A_\si^\perp$,
we see that $\imath^* f$, which is right $T_\si^\perp$-invariant since
$(\lambda, \frak t_\si^\perp)=0$, has to be a constant function. Hence $f$
is constant on that fiber, as claimed.

Since our argument is independent of the choice of fiber and element of
$V_\lambda$, we conclude that every element of $V_\lambda$ is constant on
every fiber of $\pi$. This implies (\ref{eq:tt}), and Proposition 4.3 is
now proved. \hfill $\Box$

$\spline$

We have shown that the irreducible $K$-representation with highest weight
$\lambda$ occurs in ${\cal O}({\bf L})$ if and only if 
$\lambda \in \overline{\si}$. For $\lambda \in \overline{\si}$,
it occurs with multiplicity one, and is given by
${\cal O}({\bf L})_\lambda$. 
We shall decide which of these irreducible $K$-representations are
square-integrable, in the following sense.

>From the description $G/(P,P) = (K/K_{ss}^\si) \times A_\si$, we see
that there is a $K \times A_\si$-invariant measure $\mu$ on $G/(P,P)$,
which is unique up to non-zero constant.
Given a holomorphic section $s$ of ${\bf L}$, we consider the integral
\[ \int_{G/(P,P)} <s,s> \mu \; .\]
Let $H_\omega \subset {\cal O}({\bf L})$ be the holomorphic sections
in which this integral converges. Since the Hermitian structure $<,>$
and $\mu$ are $K$-invariant, $H_\omega$ becomes a unitary $K$-representation
space. The next proposition shows which irreducible $K$-representations
occur in $H_\omega$.

Let $\lambda \in \overline{\si}$ be a dominant integral weight. Let
\[ \Phi : G/(P,P) \longrightarrow \frak k^* \]
be the moment map of the $K$-action on $(G/(P,P),\omega)$.
Recall that ${\cal O}_\lambda$ 
and ${\cal O}({\bf L})_\lambda$ 
are respectively the holomorphic functions and sections that transform by $\lambda \in \frak t_\si^*$ under the right $T_\si$-action.

$\spline$

\noindent {\bf Proposition 4.4  } {\em  Let } 
$s \in {\cal O}({\bf L})_\lambda$. {\em  Then } $s \in H_\omega$
{\em  if and only if  } $\lambda$
{\em   is in the image of the moment map. }

\noindent {\em Proof: } Let $s_o$ be the unique holomorphic section of
Proposition 4.2. Therefore, $<s_o,s_o> = e^{-F}$, where $F$ is the potential 
function of $\omega$. Since $s_o$ is non-vanishing and
$K \times T_\si$-invariant,
\[ {\cal O}({\bf L})_\lambda = {\cal O}_\lambda \otimes s_o .\]
Therefore, we are reduced to showing that $f \in {\cal O}_\lambda$ satisfies
\begin{equation}
\int_{G/(P,P)} |f|^2 e^{-F} \mu \;\; < \;\; \infty 
\label{eq:beco}
\end{equation}
if and only if $\lambda$ is in the image of $\Phi$.

Here $\mu$ is the product of a $K$-invariant measure $dk$ on
$K/K_{ss}^\si$ and a Haar measure $da$ on $A_\si$.
By the exponential map, the measure $da$ on
$A_\si$ can be identified with the Lebesgue measure $dy$ on ${\bf R}^m$,
where $m = \dim \si$.
Given $k \in K$, the left $K$-action on ${\cal O}_\lambda$,
$L_k^* : {\cal O}_\lambda \longrightarrow {\cal O}_\lambda$, is
\[ (L_k^* f)(p) = f(kp) . \]
Let $f_{1}, ... , f_{N}$ be a basis of
 ${\cal O}_\lambda$ which is orthonormal
with respect to the
(unique) $K$-invariant inner product on ${\cal O}_\lambda$. Given an element
$ f = \sum c_{i} f_{i} $
of ${\cal O}_\lambda$, 
\[ f(ky) = (L_k^*f)(y) = \sum c_{i} a_{ir} (k) f_{r} (y) ,\]
where $a_{ir} (k)$ is the $ir$th matrix coefficient of the $K$-representation 
on ${\cal O}_\lambda$ with respect to the basis above. Thus
\[ \int {|f(ky)|}^{2} dk =
\sum c_{i} \overline{c_{j}} (\int a_{ir}(k) \overline{a_{js}}(k) dk)
 f_{r}(y) \overline{f_{s} (y)} ,\]
where the integrals are taken over $K/K_{ss}^\si$.
However, by Peter-Weyl the inner integral is equal to 
\[ \frac{1}{N} {\delta}_{ij} {\delta}_{rs} , \]
(\cite{kn:chev} p.186) so the integral (\ref{eq:beco})
 reduces to 
\begin{equation}
\frac{1}{N} {\| f \|}^{2} \int_{{\bf R}^m}
 \sum {|f_{r}(y)|}^{2} e^{-F(y)} dy ,
\label{eq:fanny}
\end{equation}
where $\| f \|$ is the norm of $f$ with respect to the given 
$K$-invariant inner product 
structure on ${\cal O}_\lambda$.
 However, each of the functions $f_{r}(y)$ 
transforms under the infinitesimal
 $\frak t_\si$-action according to 
the character $\lambda \in \frak t_\si^*$, and
therefore, being holomorphic, transforms under the action of
$\frak h_\si = (\frak t_\si)_{\bf C}$
 according to the complexified character
 ${\lambda}_{\bf C} \in \frak h_\si^*$.
In particular, ${|f_{r}(y)|}^{2}$ is a constant multiple of 
$e^{2 \lambda (y)}$. Hence if $f \neq 0$, (\ref{eq:fanny})
 is a constant multiple of 
the integral 
\[ \int_{{\bf R}^m} e^{-F(y) + 2 \lambda (y)} dy .\]
However, this integral converges
 if and only if  $2 \lambda$ is
in the image of the Legendre transform of $F$ (\cite{kn:cg}, Appendix);
or equivalently if and only 
if $\lambda$ is in the image of the moment map.
This proves the proposition. \hfill $\Box$
 
$\spline$

With this result, Theorem III follows.
We see from Theorems II, III that not all irreducibles
are contained in $H_\omega$: 
The irreducible representation ${\cal O}({\bf L})_\lambda$
with highest weight $\lambda$ satisfies
${\cal O}({\bf L})_\lambda \subset H_\omega$ if and only if 
$\lambda \in \frac{1}{2} L_F(\frak a_\si) \subset \si$. This necessarily
excludes $\lambda \in \overline{\si} \backslash \si$. 
However, in the next section, we shall
see that the potential function $F$ can be constructed such that 
$\frac{1}{2} L_F(\frak a_\si) = \si$, and hence
${\cal O}({\bf L})_\lambda \subset H_\omega$ for all $\lambda \in \si$.

\newpage
\begin{center}
\section{CONSTRUCTION OF A MODEL}
\end{center}
\setcounter{equation}{0}

Let $P$ be a parabolic subgroup of $G$, and let $\si$ its corresponding
cell of dimension $m$, given in (\ref{eq:man}).
There exist dominant fundamental weights
$\alpha_1,..., \alpha_m \in \frak a_\si^*$
(\cite{kn:he2} p.498) such that 
\[ \si = \{ \sum_1^m y_i \alpha_i \;;\; y_i > 0 \} .\]
Let $F_P : \frak a_\si \longrightarrow {\bf R}$ be defined by
\begin{equation}
 F_P(v) = \sum_1^m e^{\alpha_i(v)} .
\label{eq:pot}
\end{equation}
Then $F_P \in C^\infty(\frak a_\si)$ is strictly convex, and the image of its
Legendre transform is exactly $\si$. Therefore, the moment map $\Phi$
satisfies $\Phi(A_\si) = \si$.
Extend $F_P$ to $G/(P,P)$ by $K$-invariance,
and it follows from Theorem II that
\[ \omega_P = \sqrt{-1} \partial \bar{\partial} F_P \]
is a Kaehler structure on $G/(P,P)$. 
Let ${\bf L}_P$ be the corresponding line bundle, described before.
For a dominant integral weight $\lambda$,
we let ${\cal O}({\bf L}_P)_\lambda$ denote the holomorphic sections
of ${\bf L}_P$ that transform by $\lambda$ under the right $T_\si$-action.
Let $H_{\omega_P}$ be the holomorphic sections that are 
square-integrable under (\ref{eq:int}),
so that it is a unitary $K$-representation space.
By Theorem III, 
${\cal O}({\bf L}_P)_\lambda$ is an irreducible $K$-representation
with highest
weight $\lambda$, whenever $\lambda \in \overline{\si}$. Further,
since $\Phi(A_\si) = \si$, 
${\cal O}({\bf L}_P)_\lambda \subset H_{\omega_P}$ 
whenever $\lambda \in \si$.

We repeat this geometric construction
among all the parabolic subgroups $P$ containing the fixed
Borel subgroup $B = HN$. 
In each case, we use $F_P$ in (\ref{eq:pot}) as the potential function
for the Kaehler structure $\omega_P$ on $G/(P,P)$.
Then the direct sum 
\[ \oplus_{B \subset P} \; H_{\omega_P} \]
is a model in the sense of I.M. Gelfand \cite{kn:gz}: 
a unitary $K$-representation
where all irreducibles occur with multiplicity one.

DEPARTMENT OF APPLIED MATHEMATICS, NATIONAL CHIAO TUNG UNIVERSITY, 
HSINCHU, TAIWAN.

{\em E-mail address: }
{\tt  chuah@math.nctu.edu.tw}

\end{document}